\documentclass[12pt,a4]{amsart}

\usepackage[all]{xy}
\usepackage{amscd}
\usepackage{amssymb}
\usepackage{calrsfs}

\newtheorem{theorem}{Theorem}[section]
\newtheorem{lemma}[theorem]{Lemma}
\newtheorem{corollary}[theorem]{Corollary}
\newtheorem{proposition}[theorem]{Proposition}

\theoremstyle{definition}

\theoremstyle{remark}
\newtheorem{remark}[theorem]{Remark}

\numberwithin{equation}{section}

\def\F{{\mathbb F}}

\def\Q{{\mathbb Q}}
\def\Z{{\mathbb Z}}
\def\C{{\mathbb C}}
\def\O{{\mathcal O}}

\def\p{{\mathfrak p}}
\def\q{{\mathfrak q}}
\def\U{{\mathfrak U}}
\def\L{{\mathcal L}}
\def\X{{\mathfrak X}}
\def\Y{{\mathfrak Y}}
\def\ZZ{{\mathfrak Z}}

\def\Frob{\text{\rm Frob}}

\def\Gal{\text{\rm Gal}}
\def\Hom{\text{\rm Hom}}

\def\GL{\text{\rm GL}}
\def\Spec{\text{\rm Spec}\,}
\def\Tr{\text{\rm Tr}}

\begin{document}

\title[Birational smooth minimal models]
{Birational smooth minimal models have
equal Hodge numbers in all dimensions}

\author{Tetsushi Ito}
\address{Department of Mathematical Sciences,
University of Tokyo, 3-8-1 Komaba, Meguro, Tokyo 153-8914, Japan}
\email{itote2\char`\@ms.u-tokyo.ac.jp}

\address{Max-Planck-Institut F\"ur Mathematik,
Vivatsgasse 7, D-53111 Bonn, Germany}
\email{tetsushi\char`\@mpim-bonn.mpg.de}

\thanks{The author was supported by the Japan Society for the
Promotion of Science Research Fellowships for Young
Scientists.}
\date{Sep. 20, 2002}

\begin{abstract}
This is a resume of the author's talk at
the Worhshop on Arithmetic, Geometry and Physics around
Calabi-Yau Varieties
and Mirror Symmetry (July 23-29, 2001), the Fields Institute.
The aim of this note is to prove that birational smooth
minimal models over $\C$ have equal Hodge numbers
in all dimensions by an arithmetic method.
Our method is a refinement of the method of V. Batyrev
and C.-L. Wang on
Betti numbers who used $p$-adic integration and
the Weil conjecture. Our ingredient is to use further
arithmetic results such as the Chebotarev density theorem
and $p$-adic Hodge theory.
\end{abstract}

\maketitle

\section{Introduction} 

Smooth minimal models play an important role
in birational geometry. In this paper,
we study Hodge numbers of smooth minimal models
by an arithmetic method.
Here we recall the definition of smooth minimal models.
A divisor $D$ on a smooth projective variety $X$ is called
\textit{nef} if $D \cdot C \geq 0$ for all curves $C$ in $X$.
A \textit{smooth minimal model} is
a smooth projective variety $X$
whose canonical bundle $K_X = \Omega^{\dim X}_{X}$ is nef.

The purpose of this paper is to prove the following theorem.

\begin{theorem}
\label{MainTheorem}
Let $X$ and $Y$ be smooth minimal models over $\C$.
Assume that $X$ and $Y$ are birational over $\C$.
Then, $X$ and $Y$ have equal Hodge numbers:
$$ \dim_{\C} H^i(X,\Omega_{X}^j) = \dim_{\C} H^i(Y,\Omega_{Y}^j)
   \quad \text{for all}\ i,j. $$
\end{theorem}

In this paper, we firstly compute the number of rational points
of reduction modulo $p$ of $X$ and $Y$ by $p$-adic integration.
Then we apply the following proposition
which is proved in \S \ref{Section_Galois_representation}
by combining the Weil conjecture,
the Chebotarev density theorem and $p$-adic Hodge theory
(for a variant in terms of zeta functions,
see Corollary \ref{Corollary_Lfunction}).

\begin{proposition}
\label{RationalPointsProposition}
Let $K$ be a number field.
Let $\X$ and $\Y$ be schemes of finite type over
the ring of integers $\O_K$ of $K$
whose generic fibers 
$X = \X \otimes_{\O_K} K$ and $Y = \Y \otimes_{\O_K} K$
are proper and smooth over $K$.
If $|\X(\O_K/\p)| = |\Y(\O_K/\p)|$ for all but
a finite number of maximal ideals $\p$ of $\O_K$,
then $X$ and $Y$ have equal Hodge numbers:
$$ \dim_{K} H^i(X,\Omega_{X}^j) = \dim_{K} H^i(Y,\Omega_{Y}^j)
   \quad \text{for all}\ i,j. $$
\end{proposition}

A similar statement for Betti numbers of proper smooth
varieties over a finite field is a well-known consequence
of the Weil conjecture
(see Proposition \ref{ApplicationWeilConjecture}).
We use $p$-adic Hodge theory to get information on
Hodge numbers from Galois representations.
It is likely that Proposition \ref{RationalPointsProposition}
is well-known only for specialists in arithmetic geometry.
We expect that Proposition \ref{RationalPointsProposition}
has further interesting applications in algebraic geometry.

Here we say a few words about the history of
Theorem \ref{MainTheorem}.
The case $\dim = 1$ is trivial since
birational projective smooth curves are automatically isomorphic.
The case $\dim = 2$ is also automatic because
$X$ and $Y$ are isomorphic by uniqueness of the minimal models
for surfaces.
The case $\dim = 3$ is more interesting.
By the minimal model program for threefolds, $X$ and $Y$ are
not necessarily isomorphic, but connected by a sequence of flops.
Hence we conclude, by basic properties of flops, that
$X$ and $Y$ have equal Hodge structures
(\cite{KMM},\cite{Kawamata1},\cite{Kollar}).
However, we can't do the same if $\dim \geq 4$
since the minimal model program in $\dim \geq 4$ is
still under construction.
On the other hand, Batyrev found a new way to get cohomological
properties of birational varieties in all dimensions.
He proved that birational smooth Calabi-Yau manifolds have
equal Betti numbers in all dimensions (\cite{Batyrev1}).
He used $p$-adic integration and the Weil conjecture.
Wang generalized Batyrev's result for Betti numbers
of smooth minimal models (\cite{Wang1}).

We sketch the outline of this paper.
In \S \ref{Section_minimal_models},
we prove certain geometric properties of minimal models
needed for our purpose.
In \S \ref{Section_p-adic_integration},
we recall Weil's $p$-adic integration which computes
the number of rational points of reduction modulo $p$
by integrating a gauge form on a $p$-adic manifold.
In \S \ref{Section_Galois_representation},
we recall several facts on Galois representations and prove
Proposition \ref{RationalPointsProposition}.
Finally, in \S \ref{Section_proof_main_theorem},
we prove Theorem \ref{MainTheorem} by combining above results.

\begin{remark}
After this work was completed, Fran\c cois Loeser pointed out to
the author that Theorem \ref{MainTheorem} can be obtained
by the theory of motivic integration developed
by Kontsevich \cite{Kontsevich} and Denef-Loeser \cite{DenefLoeser}
(see also \cite{Batyrev2}, \cite{Veys}, \cite{Wang1}, \cite{Wang2}).
Willem Veys kindly informed the author that
Theorem \ref{MainTheorem} is already written in
his paper \cite{Veys}, Corollary of Theorem 2.7.
Moreover, Chin-Lung Wang kindly informed the author that he
independently obtained the same result in 2000.
The author knew it in his talk at Tokyo on February 2002.
The author would like to thank them for information.
\end{remark}

\begin{remark}
Recently, the author generalized the results of this paper
and obtained an application of $p$-adic Hodge theory
to stringy Hodge numbers for singular varieties
(for details, see \cite{Ito}).
Stringy Hodge numbers were introduced by Batyrev
in \cite{Batyrev2} where he studied them
by motivic integration.
\end{remark}

\vspace{0.2in}

\noindent
\textbf{Acknowledgments.}
The author is grateful to Takeshi Saito and Kazuya Kato
for their advice and support.
He also would like to thank Takeshi Tsuji and
Shinichi Mochizuki for invaluable suggestion
on $p$-adic Hodge theory, Keiji Oguiso, Yujiro Kawamata,
and Victor Batyrev for invaluable comment and encouragement,
Yasunari Nagai for discussion about the minimal model program,
Yoichi Mieda for carefully reading a manuscript.

\section{Geometry of minimal models}
\label{Section_minimal_models}

Let $X$ and $Y$ be birational projective smooth
algebraic varieties over $\C$.
Recall that $X$ and $Y$ are called
\textit{$K$-equivalent} if there exists
a projective smooth variety $Z$
over $\C$ and proper birational morphisms
$f:Z \to X,\ g:Z \to Y$ such that
$f^{\ast} K_X = g^{\ast} K_Y$.
$$
\xymatrix{
& Z \ar[dl]_{f} \ar[dr]^{g} \\
X & & Y }
$$
Recently the notion of $K$-equivalence plays an important
role in birational geometry (\cite{Kawamata2}, \cite{Wang2}).

Although the following property of birational smooth
minimal models seems well-known for specialists in
birational geometry,
we write it for the reader's convenience
(\cite{Kollar}, \cite{Fujita}, \cite{Wang1}).

\begin{proposition}
\label{GeometricProp}
Birational smooth minimal models over $\C$ are $K$-equivalent.
\end{proposition}

\begin{proof}
Let $X$ and $Y$ be birational smooth minimal models
over $\C$.
By taking a resolution of singularities of the closure of
the graph of the birational map,
there exists a projective smooth variety $Z$ over $\C$
and proper birational morphisms $f:Z \to X,\ g:Z \to Y$.

To see $f^{\ast} K_X = g^{\ast} K_Y$,
we write the canonical bundle relation:
$$
  K_Z = f^{\ast} K_X + F + G, \qquad
  K_Z = g^{\ast} K_Y + F' + G'
$$
Here $F,G,F',G'$ are effective divisors such that
$F,F'$ are exceptional for both $f$ and $g$,
$G$ is exceptional for $f$ but not for $g$,
and $G'$ is exceptional for $g$ but not for $f$.

By symmetry, it is enough to show $F + G \geq F' + G'$.
To show this, it is enough to show $F - F' - G' \geq 0$.
Write $F - F' - G' = A - B$ such that $A,B$ are effective divisors
and they have no common component.
It is enough to show $B = 0$.
To show this, we assume $B \neq 0$.
Then we have
$$ g^{\ast} K_Y = f^{\ast} K_X + G + (F-F'-G')
                = f^{\ast} K_X + G + (A-B). $$
By taking suitable hyperplane sections
and using the Hodge index theorem,
we can take a curve $C$ in $Z$ such that $B \cdot C < 0$,
$g(C)$ is a point and $C$ is not contained in $G + A$
(\cite{Fujita}, 1.5).
Recall that $B$ is exceptional for $g$.
Since $g^{\ast} K_Y \cdot C = K_Y \cdot g(C)= 0$, we have
$$ 0 = g^{\ast} K_Y \cdot C
     = f^{\ast} K_X \cdot C + G \cdot C + A \cdot C - B \cdot C. $$
We examine the right hand side.
$f^{\ast} K_X \cdot C = K_X \cdot f(C) \geq 0$ since $K_X$ is nef.
$G \cdot C \geq 0,\ A \cdot C \geq 0$ since $C$ is not contained in
$G + A$.
Since $B \cdot C < 0$,
we conclude that the right hand side must be positive.
This is contradiction.
Hence $B = 0$ and the proof is completed.
\end{proof}

\section{$p$-adic integration}
\label{Section_p-adic_integration}

\subsection{General definitions}

Let $p$ be a prime number and $\Q_p$ be the field of
$p$-adic numbers.
Let $F$ be a finite extension of $\Q_p$,
$R \subset F$ be the ring of integers in $F$,
$m \subset R$ be the maximal ideal of $R$,
$\F_q = R/m$ be the residue field of $F$ with $q$ elements,
where $q$ is a power of $p$.
For an element $x \in F$, we define
the \textit{$p$-adic absolute value} $|x|_p$ by
$$ |x|_p = \begin{cases} q^{-v(x)} & x \neq 0 \\ 0 & x=0 \end{cases} $$
where $v : F^{\times} \to \Z$ is
the normalized discrete valuation of $F$.

Let $\X$ be a smooth scheme over $R$ of
relative dimension $n$.
We can compute the number of $\F_q$-rational points $|\X(\F_q)|$
by integrating certain $p$-adic measure on the set of
$R$-rational points $\X(R)$.
We note that $\X(R)$ is a compact and totally disconnected
topological space with respect to its $p$-adic topology.

Let $\omega \in \Gamma(\X,\Omega^n_{\X/R})$ be
a regular $n$-form on $\X$, where $\Omega^n_{\X/R}$ is
the relative canonical bundle of $\X/R$.
We shall define the $p$-adic integration of $\omega$ on $\X(R)$
as follows.
Let $s \in \X(R)$ be a $R$-rational point.
Let $U \subset \X(R)$ be a sufficiently small
$p$-adic open neighborhood of $s$
on which there exists a system of local $p$-adic
coordinates $\{ x_1,\ldots,x_n \}$.
Then $\{ x_1,\ldots,x_n \}$ defines a $p$-adic analytic map
$$ x = (x_1,\ldots,x_n) : U \longrightarrow R^n $$
which is a homeomorphism between $U$ and
a $p$-adic open set $V$ of $R^n$.
By using the above coordinates, $\omega$ can be written as
$$ \omega = f(x) \ dx_1 \wedge \cdots \wedge dx_n. $$
We can consider $f(x)$ as a $p$-adic analytic function on $V$.
Then we define the $p$-adic integration of $\omega$ on $U$
by the equation
$$ \int_{U} |\omega|_p = \int_{V} |f(x)|_p \ dx_1 \cdots dx_n, $$
where $|f(x)|_p$ is the $p$-adic absolute value of the value of $f$
at $x \in V$ and $dx_1 \cdots dx_n$ is the Haar measure on $R^n$
normalized by the condition
$$ \int_{R^n} dx_1 \cdots dx_n = 1. $$
By patching the above integration, we get
the \textit{$p$-adic integration of $\omega$ on $\X(R)$}
$$ \int_{\X(R)} |\omega|_p. $$

\subsection{$p$-adic integration of a gauge form}
\label{Subsection_gauge_form}

A \textit{gauge form} $\omega$ on $\X$ is a nowhere vanishing
global section $\omega \in \Gamma(\X,\Omega^n_{\X/R})$.
Clearly, a gauge form exists if and only if $\Omega^n_{\X/R}$
is trivial.
Therefore a gauge form exists at least Zariski locally.
The most important property of $p$-adic integration is that
the $p$-adic integration of a gauge form
computes the number of $\F_q$-rational points.

\begin{proposition}[\cite{Weil2}, 2.2.5]
\label{padicIntegration}
Let $\X$ be a smooth scheme over $R$ of relative dimension $n$
and $\omega$ be a gauge form on $\X$.
Then we have
$$ \int_{\X(R)} |\omega|_p = \frac{|\X(\F_q)|}{q^n}. $$
\end{proposition}

\begin{proof}
Let
$$ \varphi : \X(R) \longrightarrow \X(\F_q) $$
be the reduction modulo $m$ map.
For $\bar{x} \in \X(\F_q)$, $\varphi^{-1}(\bar{x})$ is
a $p$-adic open set of $\X(R)$.
Therefore, it is enough to show
$$ \int_{\varphi^{-1}(\bar{x})} |\omega|_p = \frac{1}{q^n}. $$
Let $\{ x_1,\ldots,x_n \} \subset \O_{\X,\bar{x}}$ be
a regular system of parameters at $\bar{x}$.
Then $\{ x_1,\ldots,x_n \}$ defines a system of local
$p$-adic coordinates on $\varphi^{-1}(\bar{x})$ and
$$ x = (x_1,\ldots,x_n) : \varphi^{-1}(\bar{x})
      \longrightarrow m^n \subset R^n $$
is a $p$-adic analytic homeomorphism.
Let $\omega$ be written as
$\omega = f(x) \ dx_1 \wedge \cdots \wedge dx_n$.
Since $\omega$ is a gauge form, $f(x)$ is a $p$-adic unit
for all $x \in \varphi^{-1}(\bar{x})$.
Therefore $|f(x)|_p = 1$.
Then we have
$$ \int_{\varphi^{-1}(\bar{x})} |\omega|_p
     = \int_{m^n} dx_1 \cdots dx_n = \frac{1}{q^n} $$
since $m^n$ is an index $q^n$ subgroup of $R^n$.
\end{proof}

\subsection{A slight generalization --- $p$-adic integration
of local generators of a lattice}

We shall consider Proposition \ref{padicIntegration}
if $\Omega^n_{\X/R}$ is not necessarily trivial.
By a \textit{lattice} of $\Omega^n_{\X/R}$,
we mean a locally free subsheaf $\L \subset \Omega^n_{\X/R}$
of rank 1.
We can define the $p$-adic integration of local generators of
a lattice $\L$ as follows.
If both $\L$ and $\Omega^n_{\X/R}$ are free,
take a generator $\omega$ of $\L$.
Then we have the $p$-adic integration
$$ \int_{\X(R)} |\omega|_p. $$
This value is independent of $\omega$ since
$\omega$ is unique up to multiplication by
a unit $f \in \O_{\X}^{\ast}$ and such $f$ takes
$p$-adic absolute value 1 as a $p$-adic analytic function.
Therefore we can patch them and get
the \textit{$p$-adic integration of local generators of $\L$}
which we denote by
$$ \int_{\X(R)} |\L|_p. $$
Note that if $\L = \Omega^n_{\X/R}$ and $\Omega^n_{\X/R}$ is
trivial, the above value is nothing but the $p$-adic
integration of a gauge form in \S \ref{Subsection_gauge_form}.
Therefore, we have the following
generalization of Proposition \ref{padicIntegration}

\begin{corollary}
\label{padicIntegration_lattice}
Let $\X$ be a smooth scheme over $R$ of relative dimension $n$.
Then we have
$$ \int_{\X(R)} |\Omega^n_{\X/R}|_p = \frac{|\X(\F_q)|}{q^n}. $$
\end{corollary}

\subsection{An application to $K$-equivalent varieties}

\begin{proposition}[\cite{Batyrev1}, \cite{Wang1}]
\label{CompareRationalPoints}
Let $\X,\Y,\ZZ$ be smooth schemes over $R$
of relative dimension $n$.
Assume that there exist proper birational morphisms
$\tilde{f}:\ZZ \to \X,\ \tilde{g}:\ZZ \to \Y$ such that
$\tilde{f}^{\ast} \Omega^n_{\X/R} = \tilde{g}^{\ast} \Omega^n_{\Y/R}$.
Then $|\X(\F_{q})| = |\Y(\F_{q})|$.
$$
\xymatrix{
& \ZZ \ar[dl]_{\tilde{f}} \ar[dr]^{\tilde{g}} \\
\X & & \Y }
$$
\end{proposition}

\begin{proof}
By Corollary \ref{padicIntegration_lattice}, we have
$$ \int_{\X(R)} |\Omega^n_{\X/R}|_p = \frac{|\X(\F_q)|}{q^n},
   \qquad
   \int_{\Y(R)} |\Omega^n_{\Y/R}|_p = \frac{|\Y(\F_q)|}{q^n}. $$
Since
$\tilde{f}^{\ast} \Omega^n_{\X/R} = \tilde{g}^{\ast} \Omega^n_{\Y/R}$,
we compute
\begin{align*}
\frac{|\X(\F_q)|}{q^n} &= \int_{\X(R)} |\Omega^n_{\X/R}|_p
  = \int_{\ZZ(R)} |\tilde{f}^{\ast} \Omega^n_{\X/R}|_p
  = \int_{\ZZ(R)} |\tilde{g}^{\ast} \Omega^n_{\Y/R}|_p \\
  &= \int_{\Y(R)} |\Omega^n_{\Y/R}|_p
  = \frac{|\Y(\F_q)|}{q^n}
\end{align*}
by using the change of variable formula for
$p$-adic integration twice.
Hence we have $|\X(\F_{q})| = |\Y(\F_{q})|$.
\end{proof}

\begin{remark}
By combining Proposition \ref{GeometricProp},
Proposition \ref{CompareRationalPoints},
and the Weil conjecture
(see \S \ref{Subsection_Weil_conjecture}),
we can prove the equality of Betti numbers for
birational smooth minimal models over $\C$
as in \cite{Batyrev1}, \cite{Wang1}.
However, to show the equality of Hodge numbers,
we need further arithmetic results such as the Chebotarev
density theorem and $p$-adic Hodge theory
which will be explained in
\S \ref{Section_Galois_representation}.
\end{remark}

\section{Review of Galois representations}
\label{Section_Galois_representation}

In this section, we recall several facts on
Galois representations and prove
Proposition \ref{RationalPointsProposition}.

\subsection{The Weil conjecture}
\label{Subsection_Weil_conjecture}

Let $X$ be a proper smooth variety over a finite field $\F_q$
of dimension $n$.
Fix a prime number $l$ prime to $q$.
Let $H^i_{\text{\'et}}(X_{\overline{\F}_q}, \Q_l)$ be
the $i$-th $l$-adic \'etale cohomology of
$X_{\overline{\F}_q} = X \otimes_{\F_q} \overline{\F}_q$,
where $\overline{\F}_q$ denotes an algebraic closure of $\F_q$.
Let $F : X_{\overline{\F}_q} \to X_{\overline{\F}_q}$ be
the $q$-th power Frobenius morphism.
Note that the set of fixed points of
$F$ is precisely the set of $\F_q$-rational points $X(\F_q)$.
Then, by the \textit{Lefschetz fixed point formula for \'etale
cohomology}, we have
\begin{eqnarray}
\label{Lefschetz}
|X(\F_q)| = \sum_{i=0}^{2n} (-1)^i
\ \Tr(F^{\ast};H^i_{\text{\'et}}(X_{\overline{\F}_q}, \Q_l)).
\end{eqnarray}

Moreover, by the \textit{Weil conjecture} proved by Deligne,
all eigenvalues of $F^{\ast}$ acting on 
$H^i_{\text{\'et}}(X_{\overline{\F}_q}, \Q_l)$ are
algebraic numbers and all conjugates of them
have complex absolute value $q^{i/2}$.
This is an analogue of the Riemann hypothesis
for a proper smooth variety over a finite field.

The \textit{Hasse-Weil zeta function} $Z(X,t)$
is a formal power series with coefficients in $\Q$ defined by
$$ Z(X,t) = \exp \left(
   \sum_{r=1}^{\infty} \frac{|X(\F_{q^r})|}{r} t^r \right). $$
Then, by (\ref{Lefschetz}), we have the following expression
of $Z(X,t)$
$$
Z(X,t) = \frac{P_1(X,t) \cdots P_{2n-1}(X,t)}
{P_0(X,t) P_2(X,t)\cdots P_{2n}(X,t)},
$$
where
$$ P_i(X,t) = \det(1 - F^{\ast} t;
   H^i_{\text{\rm \'et}}(X_{\overline{\F}_q}, \Q_l)). $$

Although the following application of the Weil conjecture
is well-known and weaker than our key proposition
(Proposition \ref{RationalPointsProposition}),
we note it here for reader's convenience.

\begin{proposition}
\label{ApplicationWeilConjecture}
Let $X$ and $Y$ be proper smooth varieties over $\F_q$.
If $|X(\F_{q^r})| = |Y(\F_{q^r})|$ for all $r$, then
$P_i(X,t) = P_i(Y,t)$.
In particular, by comparing the degrees, we have
$$ \dim_{\Q_l} H^i_{\text{\rm \'et}}(X_{\overline{\F}_q},\Q_l)
= \dim_{\Q_l} H^i_{\text{\rm \'et}}(Y_{\overline{\F}_q},\Q_l). $$
Therefore, if such $X$ (resp. $Y$) comes from a proper smooth
variety $\widetilde{X}$ (resp. $\widetilde{Y}$) over
a number field $K$ by modulo $p$ reduction, then
$\widetilde{X} \otimes_K \C$ and $\widetilde{Y} \otimes_K \C$
have equal Betti numbers:
$$ \dim_{\Q} H^i(\widetilde{X} \otimes_K \C, \Q)
   = \dim_{\Q} H^i(\widetilde{Y} \otimes_K \C, \Q)
   \quad \text{for all}\ i. $$
\end{proposition}

\begin{proof}
We have $Z(X,t) = Z(Y,t)$ by definition.
Hence $P_i(X,t) = P_i(Y,t)$ for all $i$ because we can recover
$P_i(X,t)$ (resp. $P_i(Y,t)$) from
$Z(X,t)$ (resp. $Z(Y,t)$) by the Weil conjecture.
The rest follows from basic properties of \'etale cohomology.
\end{proof}

\subsection{An application of the Chebotarev density theorem}

In this section, we don't need the Chebotarev density theorem
itself but need its application to $l$-adic Galois
representations.
For details, we refer Serre's book \cite{Serre}.

\begin{proposition}[\cite{Serre}, I.2.3]
\label{PropChebotarev}
Let $K$ be a number field, $m,m' \geq 1$ be integers,
and $l$ be a prime number.
Let
$$ \rho : \Gal(\overline{K}/K) \rightarrow \GL(m,\Q_l),\qquad
   \rho' : \Gal(\overline{K}/K) \rightarrow \GL(m',\Q_l) $$
be continuous $l$-adic $\Gal(\overline{K}/K)$-representations
such that $\rho$ and $\rho'$ are unramified outside a finite
set $S$ of maximal ideals of $\O_K$.
If
$$ \Tr(\rho(\Frob_{\p})) = \Tr(\rho'(\Frob_{\p}))
   \qquad \text{for all maximal ideals} \quad \p \notin S, $$
then $\rho$ and $\rho'$ have the same semisimplifications as
$\Gal(\overline{K}/K)$-representations.
Here $\Frob_{\p}$ denotes a geometric Frobenius element at $\p$
which specializes to the inverse of the $|\O_K/\p|$-th
power map on the residue field at $\p$.
\end{proposition}

\begin{proof}
We only give a sketch of the proof.
By the Chebotarev density theorem (\cite{Serre} I.2.2),
the set of conjugates of $\Frob_{\p}$ for all $\p \notin S$
is dense in
$G' = \Gal(\overline{K}/K)/(\text{Ker}\rho \cap \text{Ker}\rho')$.
Since $\rho$ and $\rho'$ are continuous representations,
$\Tr \circ \rho$ and $\Tr \circ \rho'$ are continuous maps
from $G'$ to $\Q_l$ which coincide on a dense subset of $G'$.
Hence they are equal.
Therefore, we have
$\Tr(\rho(\sigma)) = \Tr(\rho'(\sigma))$ for all
$\sigma \in \Gal(\overline{K}/K)$.
Then we have the conclusion by representation theory
(see, for example, Bourbaki, \textit{Alg\`ebre},
Ch. 8, \S 12, n${}^{\circ}$ 1, Prop 3.).
\end{proof}

\subsection{$p$-adic Hodge theory}

In this section, we recall $p$-adic Hodge theory.
Especially, we recall Hodge-Tate decomposition which is
a $p$-adic analogue of Hodge decomposition over $\C$.

Let $p$ be a prime number and $F$ be a finite extension of $\Q_p$.
Let $\C_p$ be a $p$-adic completion of an algebraic closure
$\overline{F}$ of $F$.
We define the \textit{$p$-adic Tate twists} as follows.
We define $\Q_p(0) = \Q_p,
\ \Q_p(1) = \left(\varprojlim \mu_{p^n}\right) \otimes_{\Z_p} \Q_p$,
and, for $n \geq 1$,
$\Q_p(n) = \Q_p(1)^{\otimes n},\ \Q_p(-n) = \Hom(\Q_p(n),\Q_p)$.
Moreover, we define $\C_p(n) = \C_p \otimes_{\Q_p} \Q_p(n)$,
on which $\Gal(\overline{F}/F)$ acts diagonally.
It is known that
$(\C_p)^{\Gal(\overline{F}/F)} = F$ and
$(\C_p(n))^{\Gal(\overline{F}/F)} = 0$ for $n \neq 0$
(\cite{Tate}, Theorem 2).

Let $B_{HT} = \bigoplus_{n \in \Z} \C_p(n)$
be a graded $\C_p$-module with an action of $\Gal(\overline{F}/F)$.
For a finite dimensional $\Gal(\overline{F}/F)$-representation $V$
over $\Q_p$, we define a finite dimensional graded $F$-module
$D_{HT}(V)$ by
$D_{HT}(V) = (V \otimes_{\Q_p} B_{HT})^{\Gal(\overline{F}/F)}.$
The graded module structure of $D_{HT}(V)$ is induced
from that of $B_{HT}$.
In general, it is known that
$$ \dim_{F} D_{HT}(V) \leq \dim_{\Q_p} V. $$
If the equality holds,
$V$ is called a \textit{Hodge-Tate representation}
(\cite{Tate}, \cite{Fontaine}).

\begin{theorem}[Hodge-Tate decomposition, \cite{Faltings},\cite{Tsuji}]
\label{HodgeTateConj}
Let $X$ be a proper smooth variety over $F$ and $k$ be an integer.
The $p$-adic \'etale cohomology
$H^k_{\text{\rm \'et}}(X_{\overline{F}},\Q_p)$
of $X_{\overline{F}} = X \otimes_F \overline{F}$ is
a finite dimensional $\Gal(\overline{F}/F)$-representation
over $\Q_p$.
Then, $H^k_{\text{\rm \'et}}(X_{\overline{F}},\Q_p)$
is a Hodge-Tate representation.
Moreover, there exists a canonical and functorial isomorphism
$$
 \bigoplus_{i+j = k} H^i(X,\Omega_{X}^j) \otimes_{F} \C_p(-j)
 \cong H^k_{\text{\rm \'et}}(X_{\overline{F}},\Q_p)\otimes_{\Q_p}\C_p
$$
of $\Gal(\overline{F}/F)$-representations,
where $\Gal(\overline{F}/F)$ acts on $H^i(X,\Omega_{X}^j)$
trivially and the right hand side diagonally.
\end{theorem}

For a finite dimensional $\Gal(\overline{F}/F)$-representation $V$
over $\Q_p$, we define
$$ h^{n}(V) = \dim_F(V\otimes_{\Q_p}\C_p(n))^{\Gal(\overline{F}/F)}. $$
The following lemma seems well-known for specialists in
$p$-adic Hodge theory.
But we write it here for reader's convenience.

\begin{lemma}
\label{HodgeTateLemma}
Let $W_2$ be a Hodge-Tate representation and
$$
\begin{CD}
0 @>>> W_1 @>>> W_2 @>>> W_3 @>>> 0
\end{CD}
$$
be an exact sequence of finite dimensional
$\Gal(\overline{F}/F)$-representations over $\Q_p$.
Then, $W_1$ and $W_3$ are Hodge-Tate representations
and
$$ h^n(W_2) = h^n(W_1) + h^n(W_3) = h^n(W_1 \oplus W_3) $$
for all $n$.
\end{lemma}

\begin{proof}
Since
$$
\begin{CD}
0 @>>> D_{HT}(W_1) @>>> D_{HT}(W_2) @>>> D_{HT}(W_3)
\end{CD}
$$
is exact by definition, we have
$ \dim_{F} D_{HT}(W_2)
  \leq \dim_{F} D_{HT}(W_1) + \dim_{F} D_{HT}(W_3).$
On the other hand, since $W_2$ is a Hodge-Tate representation,
we have
\begin{align*}
\dim_{F} D_{HT}(W_2) &= \dim_{\Q_p} W_2
     = \dim_{\Q_p} W_1 + \dim_{\Q_p} W_3 \\
     &\geq \dim_{F} D_{HT}(W_1) + \dim_{F} D_{HT}(W_3).
\end{align*}
Therefore, we have
$\dim_{F} D_{HT}(W_1) + \dim_{F} D_{HT}(W_3)
= \dim_{\Q_p} W_1 + \dim_{\Q_p} W_3$
and hence $W_1$ and $W_3$ are Hodge-Tate representations.
Then,
$$
\begin{CD}
0 @>>> D_{HT}(W_1) @>>> D_{HT}(W_2) @>>> D_{HT}(W_3) @>>> 0
\end{CD}
$$
is exact.
If we take the dimension of each graded quotient of
the above exact sequence, we have
$h^n(W_2) = h^n(W_1) + h^n(W_3) = h^n(W_1 \oplus W_3)$.
\end{proof}

By combining above results, we can recover the Hodge numbers
from the semisimplifications
of the $p$-adic \'etale cohomology as follows.

\begin{corollary}
\label{KeyCor}
Let $X$ be a proper smooth variety over $F$.
Then, we have
$$ \dim_{F} H^i(X,\Omega_{X}^j)
     = h^{j}(H^{i+j}(X_{\overline{F}},\Q_p)^{ss})
     \qquad \text{for all}\ i,j, $$
where $H^{i+j}(X_{\overline{F}},\Q_p)^{ss}$ denotes
the semisimplification of $H^{i+j}(X_{\overline{F}},\Q_p)$
as a $\Gal(\overline{F}/F)$-representation.
\end{corollary}

\begin{proof}
By Theorem \ref{HodgeTateConj},
if we take the dimension of the $\Gal(\overline{F}/F)$-invariant
of $H^{i+j}(X_{\overline{F}},\Q_p) \otimes_{\Q_p} \C_p(j)$,
we have
$$ \dim_{F} H^i(X,\Omega_{X}^j)
     = h^{j}(H^{i+j}(X_{\overline{F}},\Q_p)). $$
On the other hand, since $H^{i+j}(X_{\overline{F}},\Q_p)$ is
a Hodge-Tate representation,
$$ h^{j}(H^{i+j}(X_{\overline{F}},\Q_p))
   = h^{j}(H^{i+j}(X_{\overline{F}},\Q_p)^{ss}) $$
by Lemma \ref{HodgeTateLemma}.
Hence Corollary \ref{KeyCor} is proved.
\end{proof}

\begin{remark}
\label{RemarkHodgeTateConj}
A proof of Theorem \ref{HodgeTateConj} was firstly
given by Faltings (\cite{Faltings}, for recent developments
of Faltings' theory of almost \'etale extensions,
see also \cite{Faltings2}).
Tsuji gave another proof
by using de Jong's alteration (\cite{Tsuji}).
In this paper, we don't need the full version of
Theorem \ref{HodgeTateConj}.
For example,  the theorem of Fontaine-Messing is
enough for our purpose who proved Theorem \ref{HodgeTateConj}
in the case $F$ is unramified over $\Q_p$,
$\dim X < p$ and $X$ has good reduction
(\cite{FontaineMessing}).
\end{remark}

\subsection{An application to Hodge numbers}

Here we prove Proposition \ref{RationalPointsProposition}.
It is a standard consequence of the above results.

\begin{proof}[Proof of Proposition \ref{RationalPointsProposition}]
Let notation be as in Proposition \ref{RationalPointsProposition}.
Fix a prime number $l$.
Let $S$ be a sufficiently large finite set of
maximal ideals of $\O_K$ such that
$\X$ and $\Y$ are proper and smooth
over $(\Spec \O_K) \backslash S$,
$|\X(\O_K/\p)| = |\Y(\O_K/\p)|$ for all $\p \notin S$,
and $S$ contains all $\p$ dividing $l$.
Let $H^i_{\text{\'et}}(X_{\overline{K}}, \Q_l)$
(resp. $H^i_{\text{\'et}}(Y_{\overline{K}}, \Q_l)$)
be the $i$-th $l$-adic \'etale cohomology of
$X_{\overline{K}} = X \otimes_K \overline{K}$
(resp. $Y_{\overline{K}} = Y \otimes_K \overline{K}$)
on which $\Gal(\overline{K}/K)$ acts.

Let $\p$ be a maximal ideal of $\O_K$ outside $S$
and $\F_q = \O_K/\p$ be the residue field at $\p$.
Then $\X_{\overline{\F}_q} = \X \otimes_{\O_K} \overline{\F}_q$
is a proper smooth variety over $\overline{\F}_q$.
Let $F : \X_{\overline{\F}_q} \to \X_{\overline{\F}_q}$ be
the $q$-th power Frobenius morphism as in
\S \ref{Subsection_Weil_conjecture}.
Note that $\p$ doesn't divide $l$ here.
Then, by basic properties of \'etale cohomology,
$H^i_{\text{\'et}}(X_{\overline{K}}, \Q_l)$
is canonically isomorphic to
$H^i_{\text{\'et}}(\X_{\overline{\F}_q}, \Q_l)$
and the action of $F^{\ast}$ on
$H^i_{\text{\'et}}(\X_{\overline{\F}_q}, \Q_l)$
corresponds to the action of $\Frob_{\p}^{\ast}$
on $H^i_{\text{\'et}}(X_{\overline{K}}, \Q_l)$,
where $\Frob_{\p}$ is a geometric Frobenius element at $\p$
as in Proposition \ref{PropChebotarev}.
Therefore, by the Lefschetz fixed point formula for \'etale
cohomology
(\S \ref{Subsection_Weil_conjecture}, (\ref{Lefschetz})),
we have
$$
|\X(\O_K/\p)| = \sum_{i=0}^{2 \dim X} (-1)^i
\ \Tr(\Frob_{\p}^{\ast};H^i_{\text{\'et}}(X_{\overline{K}}, \Q_l)).
$$
The same is true for $\Y$.
Therefore, for $\p \notin S$, we have
$$
\sum_{i=0}^{2 \dim X} (-1)^i
\ \Tr(\Frob_{\p}^{\ast};H^i_{\text{\'et}}(X_{\overline{K}}, \Q_l))
= \sum_{j=0}^{2 \dim Y} (-1)^j
\ \Tr(\Frob_{\p}^{\ast};H^j_{\text{\'et}}(Y_{\overline{K}}, \Q_l)).
$$
Hence, by Proposition \ref{PropChebotarev},
$$
  V = \bigoplus_{i=0}^{\dim X}
  \ H^{2i}_{\text{\'et}}(X_{\overline{K}}, \Q_l)
  \ \oplus
  \bigoplus_{j=0}^{\dim Y-1}
  \ H^{2j+1}_{\text{\'et}}(Y_{\overline{K}}, \Q_l)
$$
and
$$
  W = \bigoplus_{i=0}^{\dim X-1}
  \ H^{2i+1}_{\text{\'et}}(X_{\overline{K}}, \Q_l)
  \ \oplus
  \bigoplus_{j=0}^{\dim Y}
  \ H^{2j}_{\text{\'et}}(Y_{\overline{K}}, \Q_l)
$$
have the same semisimplifications as 
$\Gal(\overline{K}/K)$-representations.

Fix a maximal ideal $\p$ of $\O_K$ outside $S$.
Let $V'$ be a simple subquotient of $V$ as
a $\Gal(\overline{K}/K)$-representation.
Then $V'$ comes from one of $H^i_{\text{\'et}}$'s.
By the Weil conjecture
(\S \ref{Subsection_Weil_conjecture}),
we can determine which cohomology group has
$V'$ as a subquotient by looking at the complex absolute
value of the eigenvalues of $\Frob_{\p}^{\ast}$ acting on $V'$.
The same is true for $W$.
Therefore, since
$H^{2 \dim X}_{\text{\'et}}(X_{\overline{K}}, \Q_l) \neq 0$
and
$H^{2 \dim Y}_{\text{\'et}}(Y_{\overline{K}}, \Q_l) \neq 0$
by Poincar\'e duality, we have $\dim X = \dim Y$.
Moreover, for each $i$, we conclude that
$H^i_{\text{\'et}}(X_{\overline{K}}, \Q_l)$
and
$H^i_{\text{\'et}}(Y_{\overline{K}}, \Q_l)$
have the same semisimplifications as 
$\Gal(\overline{K}/K)$-representations.

Now, take a maximal ideal $\q$ of $\O_K$
\textit{dividing} $l$.
Let $F$ be the completion of $K$ at $\q$.
Fix an embedding $\overline{K} \hookrightarrow \overline{F}$.
Then we have an inclusion
$\Gal(\overline{F}/F) \subset \Gal(\overline{K}/K)$.
Since $\overline{F}/\overline{K}$ is an extension
of algebraically closed fields of characteristic $0 \neq l$,
by basic properties of \'etale cohomology,
we have canonical isomorphisms
$H^i_{\text{\'et}}(X_{\overline{K}}, \Q_l)
  \cong H^i_{\text{\'et}}(X_{\overline{F}}, \Q_l)$
and
$H^i_{\text{\'et}}(Y_{\overline{K}}, \Q_l)
  \cong H^i_{\text{\'et}}(Y_{\overline{F}}, \Q_l)$
as $\Gal(\overline{F}/F)$-representations,
where $\Gal(\overline{F}/F)$ acts on
$H^i_{\text{\'et}}(X_{\overline{K}}, \Q_l)$ and
$H^i_{\text{\'et}}(Y_{\overline{K}}, \Q_l)$
by restriction.
Therefore,
$H^i_{\text{\'et}}(X_{\overline{F}}, \Q_l)$
and $H^i_{\text{\'et}}(Y_{\overline{F}}, \Q_l)$
have the same semisimplifications as 
$\Gal(\overline{F}/F)$-representations.
By Corollary \ref{KeyCor}, we finally conclude that
$X$ and $Y$ have equal Hodge numbers.
\end{proof}

\begin{remark}
If we take a sufficiently large $l$ in the above proof,
we can use the theorem of Fontaine-Messing
(Remark \ref{RemarkHodgeTateConj}, \cite{FontaineMessing}).
Therefore, we don't need the full version of
Theorem \ref{HodgeTateConj} in this paper.
\end{remark}

The following corollary is a variant of
Proposition \ref{RationalPointsProposition}
in terms of zeta functions.

\begin{corollary}
\label{Corollary_Lfunction}
Let $X$ and $Y$ be proper smooth varieties over a number
field $K$. If local zeta functions
${\zeta}_{\p}(X,s)$ and ${\zeta}_{\p}(Y,s)$ are the same for
all but a finite number of maximal ideals $\p$ of $\O_K$,
then $X$ and $Y$ have equal Hodge numbers:
$$ \dim_{\C} H^i(X,\Omega_{X}^j) = \dim_{\C} H^i(Y,\Omega_{Y}^j)
   \quad \text{for all}\ i,j. $$
\end{corollary}

\begin{proof}
Let $Z_{\p}(X,t)$ be the Hasse-Weil zeta function of
$X$ modulo $\p$ as in \S \ref{Subsection_Weil_conjecture}
for all but finitely many $\p$.
Then ${\zeta}_{\p}(X,s) = Z_{\p}(X,|\O_K/\p|^{-s})$ by definition.
The same is true for $Y$.
Therefore, Corollary \ref{Corollary_Lfunction}
is an immediate consequence of
Proposition \ref{RationalPointsProposition}.
\end{proof}

\begin{remark}
Here we note a practical difference between
Proposition \ref{RationalPointsProposition} and
Proposition \ref{ApplicationWeilConjecture}.
In Proposition \ref{ApplicationWeilConjecture},
it is sometimes possible to compute
the Hasse-Weil zeta function $Z(X,t)$
and hence Betti numbers explicitly
(for example, see \cite{Weil1}).
However, it seems very difficult to compute
the number of rational points $|\X(\O_K/\p)|$
for all but finitely many $\p$ for $X$.
Even if they are computed, there seems
no way to compute Hodge numbers explicitly.
Nevertheless, we expect that
Proposition \ref{RationalPointsProposition}
has further interesting applications
in algebraic geometry.
\end{remark}

\section{Proof of the main theorem}
\label{Section_proof_main_theorem}

In this section, we give a proof of Theorem \ref{MainTheorem}.
By Proposition \ref{GeometricProp}, it is enough to show
the following statement on Hodge numbers of $K$-equivalent
varieties.

\begin{proposition}
Let $X$ and $Y$ be birational projective smooth
algebraic varieties over $\C$.
Assume that $X$ and $Y$ are \textit{$K$-equivalent}
(see \S \ref{Section_minimal_models}),
then $X$ and $Y$ have equal Hodge numbers:
$$ \dim_{\C} H^i(X,\Omega_{X}^j) = \dim_{\C} H^i(Y,\Omega_{Y}^j)
   \quad \text{for all}\ i,j. $$
\end{proposition}

\begin{proof}
Since $X$ and $Y$ are $K$-equivalent,
there exists a projective smooth variety $Z$ over $\C$
and proper birational morphisms $f:Z \to X,\ g:Z \to Y$
such that $f^{\ast} K_X = g^{\ast} K_Y$.
$X,Y,Z,f,g$ are defined over a subfield
$K'$ of $\C$ which is finitely generated over $\Q$.
Take a number field $K$ and a variety $T$
over $K$ such that the function field $K(T)$ of $T$
is isomorphic to $K'$ over $K$.
Then, we can take varieties $X',Y',Z'$ over $K$,
morphisms $X' \to T,\ Y' \to T,\ Z' \to T$ over $K$,
and morphisms $f':Z' \to X',\ g':Z' \to Y'$ over $T$
such that the generic fibers of $X',Y',Z',f',g'$
tensored with $\C$ are $X,Y,Z,f,g$ respectively.
Then, by shrinking $T$ if necessary,
we may assume $X',Y',Z'$ are projective and smooth over $T$,
$f',g'$ are proper birational morphisms, and
${f'}^{\ast} \Omega^n_{X'/T} = {g'}^{\ast} \Omega^n_{Y'/T}$,
where $n = \dim X = \dim Y$.
Moreover, by replacing $K$ by its finite extension,
we may assume $T$ has a $K$-rational point $s \in T(K)$.
Since, for a proper smooth family of varieties
in characteristic 0, all fibers have equal Hodge numbers
(\cite{Degen}, 5.5),
we may replace $X,Y,Z,f,g$ by the fibers of
$X',Y',Z',f',g'$ at $s$.

Therefore, by changing notation,
we may assume $X,Y,Z,f,g$ are defined over
a number field $K$.
Take schemes of finite type $\X,\Y,\ZZ$ over $\O_K$
with generic fiber $X,Y,Z$.
Let $S$ be a sufficiently large finite set of
maximal ideals of $\O_K$ such that
$\X,\Y,\ZZ$ are proper and smooth over
$\U = (\Spec \O_K) \backslash S$,
and $f:Z \to X,\ g:Z \to Y$ extend to proper birational
morphisms
$\tilde{f}:\ZZ \otimes_{\O_K} \U \to \X \otimes_{\O_K} \U,
\ \tilde{g}:\ZZ \otimes_{\O_K} \U \to \Y \otimes_{\O_K} \U$
over $\U$ satisfying
${\tilde{f}}^{\ast} \Omega^n_{(\X \otimes_{\O_K} \U)/\U}
= {\tilde{g}}^{\ast} \Omega^n_{(\Y \otimes_{\O_K} \U)/\U}$.

Then for a maximal ideal $\p$ of $\O_K$ outside $S$,
the completion of $\X,\Y,\ZZ,\tilde{f},\tilde{g}$
at $\p$ satisfy the condition of
Proposition \ref{CompareRationalPoints}.
Therefore, we have
$|\X(\O_K/\p)| = |\Y(\O_K/\p)|$ for all $\p \notin S$.
By Proposition \ref{RationalPointsProposition},
we conclude that $X$ and $Y$ have equal Hodge numbers.
\end{proof}

\end{document}